\documentclass[a4paper,leqno,twoside]{amsart}%,draft
\usepackage{amsmath,amssymb,amsthm,enumerate,hyperref,color,esint}
\numberwithin{equation}{section}
\newtheorem{theorem}{Theorem}[section]

\newtheorem{proposition}[theorem]{Proposition}
\newtheorem{lemma}[theorem]{Lemma}
\newtheorem{corollary}[theorem]{Corollary}
\newtheorem{remark}[theorem]{Remark}
\newtheorem*{theorem*}{Theorem}
\theoremstyle{definition}\newtheorem{nota}{Notations}[section]

\newcommand{\Qi}{{\mathcal Q}_{\mathcal I}}
\newcommand{\Qj}{{\mathcal Q}_{\mathcal J}}
\newcommand{\Qlin}{{\mathcal Q}_{\text{\upshape lin}}}
\newcommand{\Blin}{{\mathcal B}_{\text{\upshape lin}}}
\newcommand{\Bsym}{{\mathcal B}_{\text{\upshape sym}}}
\newcommand{\mh}{\mathop{\textrm{m}_{H^-}}}
\newcommand{\mlin}{\mathop{\textrm{m}_{\lin}}}
\newcommand{\mlinrad}{\mathop{\textrm{m}_{\lin}^{\rad}}}
\newcommand{\mj}{\mathop{\textrm{m}_{\mathcal J}}}
\newcommand{\lin}{{\text{\upshape lin}}}
\newcommand{\sym}{{\text{\upshape sym}}}
\newcommand{\rad}{{\text{\upshape rad}}}

\renewcommand{\rad}{{\text{\upshape rad}}}

\def\L{{\Lambda}}
\def\l{{\lambda}}

\newcommand{\R}{\mathbb{R}}

\definecolor{darkgreen}{rgb}{0.0, 0.5, 0.2} 
\definecolor{purple}{rgb}{0.5, 0.0, 0.5}
\newcommand{\AL}{\color{purple}}

\newcommand{\jcal}{{\mathcal J}}

\newcommand{\remove}[1]{}

\def\sideremark#1{\ifvmode\leavevmode\fi\vadjust{\vbox to0pt{\vss% the remark
 \hbox to 0pt{\hskip\hsize\hskip1em%                          will appear only
 \vbox{\hsize2.1cm\tiny\raggedright\pretolerance10000%          on the side
  \noindent #1\hfill}\hss}\vbox to15pt{\vfil}\vss}}}%

	\title{Some remarks about the Morse index for convex Hamiltonian systems}
	\author[A.~L.~Amadori]{Anna Lisa Amadori}
	\thanks{\small{This work was partially supported by Gruppo Nazionale per l'Analisi Matematica, la Probabilit\`a e le loro Applicazioni (GNAMPA) of the Istituto Nazionale di Alta Matematica (INdAM) and by the\textquotedblleft Geometric-Analytic Methods
			for PDEs and Applications\textquotedblright (GAMPA)  project CUP I53D23002420006- funded by European Union - Next Generation EU within the PRIN 2022 program (D.D. 104 - 02/02/2022 Ministero dell’Universit\`a e della Ricerca). This manuscript reflects only the author’s views and opinions and the Ministry cannot be considered responsible for them.}}
	\date{\today}
	\address{\small{Dipartimento di Scienze e Tecnologie, Universit\`a di Napoli ``Parthenope", Centro Direzionale di Napoli, Isola C4, 80143 Napoli, Italy. \newline \texttt{annalisa.amadori@uniparthenope.it}}}
\begin{document}
\begin{abstract}
	\noindent
	We investigate the (linearized) Morse index of solutions to  Hamiltonan systems, with a focus on convex Hamiltonians functions and sign-changing radial solutions.
	For strongly coupled systems, we describe the profile of the radial solutions and give an estimate of their Morse index. 
\end{abstract}
	
\maketitle

\section{Introduction}
This paper deals with elliptic systems of Hamiltonian type
\begin{equation}\label{Ham}
	\begin{cases}	-\Delta u = H_v(x,u,v)& \text{ in } \Omega , \\	
		-\Delta v = H_u  (x,u, v) \qquad & \text{ in } \Omega ,\\
		u=v=0 & \text{ on } \partial\Omega,
	\end{cases}
\end{equation}
where $\Omega$ is a bounded smooth domain in $\R^N$, $N\ge 2$, and $H$ is continuous with respect to $x\in \Omega$ and  $C^2$ with respect to $(u,v)\in \R^2$. We are particularly interested in systems where the Hamiltonian function $H$ is convex w.r.t. the variables $u, v$, that is
\begin{align} \label{H1} H_{uu}  , \,  H_{vv}   \ge 0 \  \text{ and }   \ H_{uu} H_{vv}  \ge \left( H_{uv}\right)^2 \end{align}on $\Omega\times \R^2$. 
Notice that the convexity of $H$ implies that the system is cooperative. \\
 The system \eqref{Ham} is the Euler-Lagrange equation of the action functional 
\begin{align}\label{energy} 
	\mathcal{I}(u,v) = &\int_{\Omega} \nabla u \nabla v dx  -\int_{\Omega} H(x,u,v) dx ,  %=\frac{1}{2} \mathcal Q_o (u,v) -\int_{\Omega} H(u,v) dx ,
\end{align}
meaning that weak solutions are critical points for $\mathcal I$.
The principal part 
\[ \mathcal S_o (u,v) = \int_{\Omega} \nabla u \nabla v dx \]
is strongly indefinite, as $\left(H^1_0\right)^2$ splits into the direct sum of the two sets  \[H^{\pm} := \left\{ (\phi, \pm \phi) \, : \, \phi \in H^1_0\right\}\]
%\begin{align}\label{def:Hpm}H^{\pm} := \left\{ (\phi, \pm \phi) \, : \, \phi \in H^1_0\right\}\end{align}
where $\pm \mathcal S_o$ is coercive, respectively.
Various approaches have been proposed to get rid of this difficulty and establish existence of solutions in  suitable functional and variational frameworks, together with qualitative and symmetry properties. We refer to the surveys \cite{DeFig},  \cite{BDT} and the references therein 
for a comprehensive treatment of this topic. We also mention \cite{BDRT} concerning least energy sign changing solutions. To our purpose we emphasize that every choice of the  functional setting produces continuous solutions and justifies the following assumption
\begin{equation} \label{Huv} \begin{split}
H_{uu}(x,u(x), v(x))  , \  H_{vv} (x,u(x), v(x)) ,  \  H_{uv}(x,u(x), v(x)) \\ \mbox{ are bounded on $\Omega$.} %\in L^{\infty}(\Omega). 
\end{split}\end{equation}

Our focus is on the Morse index of solutions. A first remark is that the quadratic form associated to the second derivative of $\mathcal I$ is
\begin{align}\label{def_Qi} 	\Qi(\phi,\psi) =	
	\langle\mathcal{I}''(u,v)(\phi,\psi), (\phi,\psi)\rangle  \\ \nonumber  =  2\int_{\Omega} \nabla \phi \nabla \psi dx  -  \int_{\Omega} D^2 H(\phi,\psi) \!\cdot\! (\phi,\psi)  dx  \le 2\int_{\Omega} \nabla \phi \nabla \psi dx .
\end{align}
by the convexity of $H$ % \eqref{H1}.
(here and henceforth we write $D^2 H$ meaning the Hessian matrix of $H(x,\cdot)$ with respect to the variables $(u,v)\in \R^2$).
So, $\Qi$ is negative defined on the half-space $H^-$ and the \textquotedblleft natural\textquotedblright\ Morse index of any weak solution, meaning the maximal dimension of a subspace where $-\mathcal I''$ is coercive, is not finite.
Several notions have been proposed in the literature to overcome this problem and define a meaningful finite index. Among others, we mention the relative Morse index proposed by Abbondandolo \cite{Abb}, which measures the relative dimension of the negative eigenspace of $\mathcal I''$ with respect to $H^-$, and the reduced Morse index, which computes the standard Morse index of a reduced functional, corresponding to a scalar equation. The reduced functional and Morse index have been used, for instance, to estimate the relative Morse index from below and obtain a Liouville type result  \cite{Ram09}, and to construct multiple spike solutions \cite{Ram09}.
A different approach consists in looking at the linearized Morse index, which is not directly related to the action functional $\mathcal I$, but  provides a finite index which can be used, for instance, to prove that least energy solutions are foliated  Schwartz  symmetric, see \cite{DP13SiamJMA}.

In this paper we point out some general properties of the linearized Morse index for Hamiltonian systems, with an emphasis on the convex case \eqref{H1}.
In Section 2, after recalling the definitions and basic properties, we show that the computation of the linearized Morse index reduces to the investigation of  some scalar eigenvalue problems.
\begin{proposition}\label{prop:solosym+}
	Let $(u,v)$ be a solution to \eqref{Ham} satisfying \eqref{Huv}.
	Its linearized Morse index is the sum of the number (counting multiplicity) of negative eigenvalues of two scalar problems
	\begin{equation}\label{PA_sym+intro}\begin{cases}
			-\Delta \phi + a \phi = \mu \phi \quad & \text{ in } \Omega ,  \\
			\phi=0  & \text{ on } \partial \Omega ,   
	\end{cases}\end{equation}
	\begin{equation}\label{PA_sym-intro}\begin{cases}
-\Delta \psi + b\psi= \nu \psi \quad & \text{ in } \Omega ,  \\
\psi=0  & \text{ on } \partial \Omega ,   
\end{cases}
\end{equation}
where $\displaystyle  a(x):=-\frac{1}{2}\big(H_{uu}(x,u(x),v(x)) + H_{vv}(x,u(x),v(x))  \big) - H_{uv}(x,u(x),v(x))$, \\
and $\displaystyle  b(x):=\frac{1}{2}\big(H_{uu}(x,u(x),v(x)) + H_{vv}(x,u(x),v(x))  \big) - H_{uv}(x,u(x),v(x))$. \\
Moreover the negative eigenspace is spanned by functions of type $(\phi,\phi)$ where $\phi$ solves \eqref{PA_sym+intro}, and $(\psi,-\psi)$ where $\psi$ solves \eqref{PA_sym-intro}.

Under the convexity assumption \eqref{H1},  the linearized Morse index is equal to the number (counting multiplicity) of negative eigenvalues of \eqref{PA_sym+intro}, and the negative eigenspace is contained in $H^+$.
\end{proposition}

 Next, in Section 3, we focus on radial solutions and give a formula for their Morse index which makes use of the notion of singular eigenvalues and extends to Hamiltonian systems the analogous formula for scalar equations, see \cite{AGNARWA}.
 For convex problems, we obtain the following
 
 \begin{proposition}\label{prop:morse-index-formula}
 	Consider the system \eqref{Ham} with $\Omega$ radially symmetric  and $H=H(|x|, u, v)$ satisfying \eqref{H1}. If $(u,v)$ is a radial solution satisfying \eqref{Huv}, then its linearized Morse index is given by
 \begin{equation}\label{quattro} \mlin(u,v) = {\mlinrad} +\sum\limits_{k =1}^{\mlinrad} \sum\limits_{j=0}^{M_k} N_j, \end{equation}
 where $\mlinrad$ stands for the radial linearized Morse index,  $N_j=\frac{(N+2j-2)(N+j-3)!}{(N-2)!j!}$,  \\  $M_k= \min\left\{ n\in \mathbb N  : n\ge  \sqrt{\left(\frac{N-2}{2}\right)^2-\widehat\L^{\rad}_k}-\frac{N}{2} \right\}$, \\
  and $\widehat\L^{\rad}_k$  are the singular radial eigenvalues, characterized by
 \begin{equation*}\label{PA_sym+_sing_rad_intro}
 \begin{cases} 
 -\left( r^{N-1} \xi'\right)'- r^{N-1} a \xi= r^{N-3} \widehat\L^{\rad}\xi  , \\
 \xi \in \widehat{H}^1_{0,\rad} .\end{cases}
 \end{equation*}
\end{proposition}
 
In Section 4, we deal with a  class of strongly coupled Hamiltonian systems with $H=H(u,v)$ and
\begin{align}\label{A}\tag{H1}
H_u/u > 0 &  \text{ and }  H_v /v > 0  , \\
\label{B}\tag{H2}
H_{uu} > \frac{H_u - H_{uv} v}{u}  > 0 &  \text{ and }  H_{vv}> \frac{H_v - H_{uv} u}{v} > 0   , 
\end{align}
for every $u , v \neq 0$. % {\textcolor{red}{and almost every $x\in \Omega$}}. 
These assumptions are satisfied, for instance, by the Lane-Emden system
\begin{equation*}\label{LE}
\begin{cases}	-\Delta u = |v|^{q-1} v& \text{ in } \Omega , \\	
-\Delta v = |u|^{p-1}u \qquad & \text{ in } \Omega ,\\
u=v=0 & \text{ on } \partial\Omega,
\end{cases}
\end{equation*}
for $p , q >1$ and, more generally, by every fully coupled system for which \[H(u,v)= F(v)+G(u)\] with $F''(v) > F'(v)/v  >0$ and $G''(u) > G'(u)/u > 0$  for every $u, v \neq 0$.

We prove the following estimate. 
\begin{theorem}\label{teo-stima-Morse}
Consider the system \eqref{Ham} with $\Omega$ radially symmetric %$=\{ x\in \R^N : |x|<R\}$ 
and $H=H(u, v)$ satisfying		\eqref{A} and \eqref{B}. If $(u,v)$ is a classical radial solution  and $m$ is the number of nodal zones of $u$ and $v$, then 
	\begin{equation}\label{uno} \mlinrad \ge m \end{equation}
	and 
\begin{equation}\label{due}\mlin(u,v) \ge \mlinrad +(m-1) N
	\ge m+(m -1) N . \end{equation}
\end{theorem}

Theorem \ref{teo-stima-Morse} is a consequence of the just stated general properties, but also of the strong coupling between $u$ and $v$ yielded by hypothesis \eqref{A}.
We think that this fact can be of some interest for itself, therefore we give right now a specific statement.
%Let $\Omega$ be either a ball $\{x\in\R^N : |x|<R\}$  or a spherical shell $\{x\in\R^N : \delta<|x|<R\}$, and $(u,v)$ a radial function. If $m$ and $n$ are the numbers of nodal zones of $u$ and $v$, respectively, we write $s_0=t_0 = 0$ or $s_0=t_0=\delta$ according if $\Omega$ is a ball or a spherical cell, $s_m=t_n=R$, $s_1<\dots s_{m-1}$ and $t_1<\dots t_{n-1}$ for the internal zero points  of $u(|x|)$ and $v(|x|)$,  and $\sigma_0<\sigma_1<\dots \sigma_{m-1}$ or, respectively, $\tau_0< \tau_1<\dots \tau_{n-1}$ for the respective maximum/minimum points in each nodal zone.

\begin{proposition}\label{profilo} Let $\Omega$ be radially symmetric, $H=H(|x|, u, v)$ continuous w.r.t. $x$ and satisfying \eqref{A} for almost every $x\in \Omega$, and $(u,v)$ a classical radial solution to \eqref{Ham}. Then the following holds
	\begin{enumerate}[A.]
		\item If one between $u$ and $v$ is strictly positive (or negative), then both $u$ and $v$ are strictly positive (or negative), and have only one critical point, which coincides with the origin if $\Omega$ is a ball.
		\item If one between $u$ and $v$ has exactly $m\ge 2$ nodal zones and $m-1$ internal zeros,  then  both $u$ and $v$ have exactly $m$ nodal zones and $m-1$ internal zeros, they  have the same sign in their first nodal zone (and therefore in the following ones), and have exactly one critical point inside each nodal zone. 
	Moreover the nodal zones are intertwined, meaning that the $i^{th}$ nodal zones of $u$ and $v$ have a non-empty intersection, which contains both the critical points of $u$ and $v$.
\end{enumerate}
%	\begin{align}\label{punti-critici}	\sigma_i, \tau_i  \in ( s_i, s_{i+1}) \cap ( t_i, t_{i+1} ) \qquad \text{ for } i=0,\dots m-1	\end{align}	and, when $m\ge 3$, also	\begin{align}	\label{zone-nodali}	s_{i-1} < t_i < s_{i+1} , \qquad t_{i-1} < s_i < t_{i+1} \qquad \text{ for } i=1,\dots m-1 .	\end{align}
\end{proposition}

\section{The linearized Morse index and the symmetric eigenvalues}
The notion of linearized Morse index is focused on the linearization of problem \eqref{Ham} near a solution $(u,v)$, i.e.
\begin{equation}\label{Ham-lin}
\begin{cases}	-\Delta \phi = H_{uv} \phi + H_{vv}\psi& \text{ in } \Omega , \\	
-\Delta \psi = H_{uu}\phi + H_{uv}\psi \qquad & \text{ in } \Omega ,\\
\phi=\psi=0 & \text{ on } \partial\Omega,
\end{cases}
\end{equation}
and on the related bilinear and quadratic form
\begin{align}\label{def_bil_lin}
\Blin\left((\phi,\psi), (\xi,\eta)\right) =  \int_{\Omega} \left(\nabla\phi\nabla\xi + \nabla \psi\nabla\eta\right) dx %\\ \nonumber 
- \int_{\Omega} \left(  H_{uu} \phi \eta +H_{vv} \psi\xi +H_{uv}\left(\phi \xi +\psi \eta \right) \right)dx,
\\ \label{def_quad_lin}
\Qlin(\phi,\psi) =  \int_{\Omega} \left(|\nabla\phi|^2 + |\nabla \psi|^2\right) dx - \int_{\Omega} \left( (H_{uu}+H_{vv}) \phi \psi  +  H_{uv} ( \phi^2+\psi^2 \right) dx,
\end{align}
defined respectively on $\left(H^1_0\right)^2\times  \left(H^1_0\right)^2$ and $\left(H^1_0\right)^2$.

The linearized Morse index is the maximal dimension of a subspace of 
$\left(H^1_0\right)^2$ where the quadratic form $\Qlin$ is negative defined. We will denote it by $\mlin(u,v)$, henceforth.

Let us remark that $\Qlin$ differs from the related quadratic form $\Qi$ introduced in \eqref{def_Qi}.
Interestingly, 
\begin{equation}\label{Qlin=-Qi}\Qlin(\phi,-\phi) = 2\int_{\Omega} |\nabla\phi|^2 dx +   \int_{\Omega} \left(H_{uu}+H_{vv} - 2H_{uv} \right) \phi^2 dx
= - \Qi(\phi,-\phi) ,\end{equation}
while in the orthogonal complement $H^{+}$ the two quadratic forms coincide:
\begin{equation}\label{Qlin=Qi}\Qlin(\phi,\phi) = 2\int_{\Omega} |\nabla\phi|^2 dx -   \int_{\Omega} \left(H_{uu}+H_{vv}+ 2H_{uv} \right) \phi^2 dx
=  \Qi(\phi,\phi) .\end{equation}

Hence the notion of linearized Morse index can provide a finite index which is meaningful also from the point of view of the action functional.
In order to take advantage of the theory of bounded self-adjoint operators, however, one can not rely on the bilinear form $\Blin$, which is not symmetric, but instead on its symetrization
\begin{align}\label{def_bil_sim} \begin{split}
\Bsym\left((\phi,\psi), (\xi,\eta)\right) = \int_{\Omega} \left(\nabla\phi\nabla\xi + \nabla \psi\nabla\eta\right) dx \\ - \frac{1}{2}\int_{\Omega} \left(\Delta H (\phi\eta+ \psi\xi) +  2H_{uv} ( \phi\xi+\psi\eta )
\right) dx,\end{split}
\end{align} 
which provides the same quadratic form as
\[\Qlin(\phi,\psi) =\Bsym\left((\phi,\psi), (\phi,\psi)\right) .\]
The standard spectral decomposition theory applies 
provided that there exists $k\in \R$ such that the continuous, symmetric bilinear form $\Bsym+ k \langle, \rangle_{(L^2)^2}$ is coercive. 
It certainly holds true if the solution $(u,v)$ fulfills \eqref{Huv}.
In that way there exists a non-decreasing sequence of eigenvalues $\Lambda_n^{\sym}\to +\infty$.
The first eigenvalue is the  minimum of the Rayleigh quotient
\begin{align}\label{Reyleigh1}
{\Lambda}^{\sym}_1:= \inf \left\{ \dfrac{\Qlin(\phi,\psi)}{\|\phi\|_2^2+\|\psi\|_2^2}\, : \, \phi,\psi \in H^1_0 \right\}
\end{align}
which is attained by a nontrivial function $(\phi_1,\psi_1)$ (an eigenfunction) which satisfies
\[	\Bsym \left((\phi_1,\psi_1), (\xi,\eta)\right) =\Lambda_1^{\sym} \int_{\Omega} \left(\phi_1\xi+\psi_1\eta\right) dx
\]
for every $\xi, \eta \in H^1_0$, i.e. solves in weak sense
\begin{equation}\label{PA_sym}
\begin{cases} 
-\Delta \phi - H_{uv}  \phi- 
\frac{1}{2}\Delta H \psi= \L \phi \qquad & \text{ in } \Omega  ,  \\
-\Delta  \psi -\frac{1}{2}\Delta H  \phi - H_{uv} \psi= \L \psi& \text{ in } \Omega , 
\\
\phi=\psi=0 & \text{ on } \partial \Omega \end{cases}
\end{equation}
for $\L=\L^{\sym}_1$. Conversely, any weak solution of \eqref{PA_sym} with $\L=\L^{\sym}_1$ realizes the minimum in \eqref{Reyleigh1}.
Next, the following eigenvalues are defined by 
\begin{align*}
{\Lambda}^{\sym}_{k}:= \inf \left\{ \dfrac{\Qlin(\phi,\psi)}{\|\phi\|_2^2+\|\psi\|_2^2}\, : \, \phi,\psi \in H^1_0 \, \int_{\Omega} \left(\phi\phi_n +\psi\psi_n \right) dx = 0 \text{ if } n=1,\dots k-1 \right\}
\end{align*}
with $k\ge 2$, and are attained by eigenfunctions $(\phi_k,\psi_k)$ which solve \eqref{PA_sym} for $\L=\L^{\sym}_k$.
Eventually, the family of eigenfunctions forms an orthogonal basis for $\left(H^1_0\right)^2$ and induces a splitting
\[ \left(H^1_0\right)^2 = W^-\oplus W^+\oplus W^0 %{\mathrm{Ker}}(\mathcal S_{\text{sym}})
.\]
Here $W^{+}$, $W^-$ and $W^0$ are the spaces spanned by the eigenfunctions $(\phi_k,\psi_k)$ with $\Lambda_k^{\text{sym}} <0$ (resp., $>0$, $=0$), so that  $\Qlin$ is negative (resp., positive) defined on $W^{-}$ ( resp, $W^+$).
Moreover $W^-$ has finite dimension because $\Lambda_k^{\sym}\to +\infty$.
Summing up, the linearized Morse index of a solution $(u,v)$ is the dimension of $W^-$, i.e. the number of negative eigenvalues of problem \eqref{PA_sym}, each counted with multiplicity.
We refer to \cite{DP13SiamJMA} for rigorous proofs and more details.

\subsection{Reduction to a scalar eigenvalue problem}
Projecting the eigenvalue problem \eqref{PA_sym} into $H^+\oplus H^-$ gives rise to a decoupled system.
If $(u,v)$ is a solution to \eqref{Ham}, we introduce the functions
\begin{align}\label{def_a}
a(x)=&-\frac{1}{2} \left(H_{uu}(x,u(x),v(x)) + H_{vv}(x,u(x),v(x))\right) - H_{uv}(x,u(x),v(x)) , \\
\label{def_b}
b(x)=&\frac{1}{2}\big(H_{uu}(x,u(x),v(x)) + H_{vv}(x,u(x),v(x))  \big) - H_{uv}(x,u(x),v(x)) ,
\end{align}
and the scalar eigenvalue problems
	\begin{equation}\label{PA_sym+}\begin{cases}
-\Delta \xi + a \xi = \mu \xi \quad & \text{ in } \Omega ,  \\
\xi=0  & \text{ on } \partial \Omega ,   
\end{cases}\end{equation}
\begin{equation}\label{PA_sym-}\begin{cases}
-\Delta \eta+ b\eta= \nu \eta \quad & \text{ in } \Omega ,  \\
\eta=0  & \text{ on } \partial \Omega ,   
\end{cases}
\end{equation}

It is very easy to prove the following interesting property.
\begin{lemma}\label{lem:PAdiag}
	Let $(u,v)$ a weak solution to \eqref{Ham} satisfying \eqref{Huv}. %, with any function  $H\in C\left(\Omega, C^2(\R^2)\right)$.
	A number $\Lambda$ is a symmetric eigenvalue if and only if it is an eigenvalue for at least one between \eqref{PA_sym+} or \eqref{PA_sym-}, and the related eigenfunctions of \eqref{PA_sym} are linear combinations of functions of type $(\xi ,\xi)$ and/or $(\eta ,-\eta)$, respectively.
	In particular the multiplicity of $\Lambda$ as an eigenvalue of \eqref{PA_sym} is the sum  of its  multiplicity according to \eqref{PA_sym+} and  \eqref{PA_sym-}, counting $0$ if $\Lambda$ is not an eigenvalue for \eqref{PA_sym+}, or respectively \eqref{PA_sym-}.
\end{lemma}
\begin{proof}
For $(\phi,\psi)\in \left( H^1_0\right)^2$ we write $\xi= \frac{1}{2} (\phi+\psi)$ and   $\eta =\frac{1}2(\phi-\eta)$, 
so that
$P_{H^+}(\phi,\psi)=(\xi,\xi)$,  $P_{H^-}(\phi,\psi)=(\eta,-\eta) $,
and $(\phi,\psi)= (\xi+\eta, \xi-\eta)$.
By a trivial computation, the system \eqref{PA_sym} translates into the decoupled system
\begin{equation*}\label{epHdiag}\begin{cases}
-\Delta \xi + a \xi= \Lambda \xi \qquad &\text{ in } \Omega  , \\
-\Delta \eta +b \eta = \Lambda \eta & \text{ in } \Omega , 
\\
\xi=\eta=0  & \text{ on } \partial \Omega 
\end{cases}
\end{equation*}
and the claim follows.
\end{proof}

%In particular $\Lambda_1 = \min\{ \mu_{1},\nu_{1}\}$. If the system \eqref{H} is fully coupled, (i.e. the assumptions \eqref{H-coop} and \eqref{H-compl-acc} hold true), $a\le b$ with strict inequality on a non-null set, therefore $\Lambda_1=\mu_1<\nu_1$ and $W_1=(\phi_1,\phi_1)$. (As known, the first eigenvalue is simple and the related eigenfunction has positive component). 
%But  it is not true that all the other eigenfunctions change sign, because there exists $n$ such that $\Lambda_n=\nu_1$ and its eigenfunction has fixed sign component, but with opposite sign.
%\\
%\oss{ In generale  $\Lambda_n \le \min\{\mu_n , \nu_n\}$, dunque se non si verifica mai che $\mu_k=\nu_h$ possiamo ancora dire che l'$n$-esima autofunzione ha al pi\`u $n$ zone nodali, altrimenti andando a fare le combinazioni lineari fra $\phi_k$ e $\psi_h$ non so cosa pu\`o succedere.
%Mi pare che, con un po' di regolarit\`a, si pu\`o dimostrare  che se il sistema \`e fortemente accoppiato allora $\mu_k\neq\nu_h$ (perch\'e $a<b$ su un insieme di misura positiva). }\medskip

\begin{proof}[Proof of Proposition \ref{prop:solosym+}]
The first statement, concerning any Hamiltonian function, follows readily by Lemma \ref{lem:PAdiag} and the characterization of the linearized Morse index by means of the symmetric eigenvalues.

 It remains to check that  the convexity assumption \eqref{H1} implies that
$b\ge 0$ pointwise, so that the eigenvalue problem \eqref{PA_sym-} only has nonnegative eigenvalues. \\
If $ H_{uu}=0$, then $H_{uv}=0$ and $b=H_{vv} \ge 0$. Otherwise  $H_{uu}> 0$, and
\begin{align*}
2 H_{uu} b = \left(H_{uu}\right)^2 + H_{uu} H_{vv} - 2H_{uu} H_{uv}  \underset{\eqref{H1}}\ge 	  \left(H_{uu} - H_{uv}\right)^2 
\ge 0 .
\end{align*}
\end{proof}

In particular in the convex case,  Proposition \ref{prop:solosym+}  shows that the negative eigenspace of the quadratic form $\mathcal Q_{\lin}$, i.e. $W^-$, is contained in $H^+$. Moreover, as far as negative eigenvalues are concerned, the symmetric eigenfunctions inherit all the well-known properties of the eigenfunctions for compact scalar problems: the first eigenvalue is simple and has a positive eigenfunction, and eigenfunctions related to the $n^{th}$ eigenvalues have at most $n$ nodal zones (componentwise).

\subsection{Linearized and relative Morse index}
Here we briefly recall the notion of relative Morse index introduced by Abbondandolo \cite{Abb} for Hamiltonian systems, and compare it with the linearized Morse index.

If $(u,v) \in \left(H^1_0\right)^2$ is a weak solution of \eqref{Ham}, 
we  write $\mathcal S$ for the self-adjoint realization of the symmetric bilinear form defined by  $\mathcal I''(u,v)$ on $H^1_0$, i.e. 
\begin{align*}
	\langle\mathcal S (\phi,\psi), (\xi,\eta)\rangle = \langle(\phi, \psi) , \mathcal S (\xi, \eta)\rangle=
	\langle\mathcal I''(u,v)(\phi,\psi), (\xi,\eta)\rangle \\
	= \int_{\Omega}\left(\nabla\phi \nabla\eta +\nabla\psi\nabla\xi \right) dx -   \int_{\Omega} D^2H(u,v) (\phi,\psi)\cdot (\xi,\eta)  dx .
\end{align*}
$\mathcal S$ is the difference between the strongly indefinite, but invertible, operator 
\begin{align*} %\mathcal B_o\left( (u,v), (w,z) \right) =  
	\langle\mathcal S_o(\phi,\psi) , (w,z)\rangle
	= & \int_{\Omega} \left(\nabla \phi \nabla z +  \nabla \psi \nabla w \right) dx , \intertext{and the operator}
	\langle\mathcal S_{H}(\phi,\psi) , (w,z)\rangle =& \int_{\Omega} D^2H(u,v) (\phi,\psi)\!\cdot\!(w,z) dx  ,
\end{align*}
which is compact under suitable  assumptions on the Hamiltonian $H$ (see, for instance, \cite[Proposition 3.2.1]{Abb}).
So $\mathcal S$ is a Fredholm operator and  determines an unique $\mathcal S$-invariant orthogonal splitting 
\[ \left(H^1_0\right)^2 = V^- \oplus V^+ \oplus {\mathrm{Ker}}(\mathcal S) \]
in such a way that $\Qi$ is negative (resp., positive) defined on $V^-$ (resp., $V^+$). 
As mentioned in the Introduction,  the dimension of $V^-$, i.e. the \textquotedblleft natural\textquotedblright\ Morse index, is infinite for every solution $(u,v)$ under assumption \eqref{H1}.

The \emph{relative Morse index} is the relative dimension of $V^-$ with repect to $H^-$, that is 
\begin{equation}\label{def_rel_morse} 
	{\mathrm{m}}_{H^-}(u,v):={\dim }(V^-, H^-):={\dim }(V^-\cap H^+) - {\dim }((V^-)^{\perp}\cap H^-) .
\end{equation}
%Both the dimensions appearing in this definition are finite because they are kernels of the compact operators $P_{V^-} P_{H^+}$ and $P_{H^-} P_{(V^-)^{\perp}} $, where  $P$ stands for orthogonal projection.
The convexity assumption \eqref{H1} implies that $(V^-)^{\perp}\cap H^-=\{(0,0)\}$, then in this case \eqref{def_rel_morse} simplifies into
\begin{equation}\label{rel_morse} {\mathrm{m}}_{H^-}(u,v)={\dim }(V^-\cap H^+).
\end{equation}
In any case
\begin{equation}\label{rel_morse<} {\mathrm m}_{H^-} (u,v) \le {\mathrm m}_{\lin} (u,v) .\end{equation}
Indeed by definition  ${\mathrm m}_{H^-} (u,v) \le {\dim }(V^-\cap H^+)$, and obviously $Q_{\mathcal I}$ is negative on $V^-\cap H^+$. By \eqref{Qlin=Qi}, also $ \Qlin $ is negative, hence ${\dim }(V^-\cap H^+)\le \mlin (u,v)$.
%{\textcolor{red}{	Moreover when the  Hamiltonian function is convex, then $\mlin(u,v) = \mj(u,v)$ if and only if 	\begin{equation}\label{iff}		\sup \left\{ \frac{\Qi(\xi+\eta, \xi-\eta)}{\|\xi\|_2^2 + \|\eta\|_{2}^2} \, : \, \xi\in W^-\setminus\{0\} ,\, \eta \in H^1_0 \setminus\{0\} \right\} < 0 .		\end{equation}	secondo me per dire questo bisogna passare dall'indice ridotto}}

We explicitly remark that, although in the convex case $\Qi$ is negative on the subspaces $W^-\subset H^+$ and $H^-$, one can not infer that $\Qi$ is negative on $W^-\oplus H^-$.
Indeed 	$\Qlin(\xi+\eta, \xi-\eta) = \Qlin(\xi,\xi)+ \Qlin(\eta, -\eta)$, {but}
\begin{align*}
 	\Qi(\xi+\eta, \xi-\eta) & = \Qi(\xi,\xi)+ \Qi(\eta, -\eta)+2\langle {\mathcal I}''(\xi,\xi), (\eta,-\eta)\rangle
	\\
	& =\Qi(\xi,\xi)+ \Qi(\eta, -\eta)+ 2 \int_\Omega \left(H_{vv}-H_{uu}\right) \xi\eta dx.
\end{align*}
%{\taglia COSA MANCA PER L'UGUAGLIANZA
%\begin{proposition}\label{1.2}	Under assumption \eqref{H1}, for every weak solution of \eqref{Ham} we haven\[ {\mathrm m}_{H^-} (u,v) = \mlin (u,v) \]\end{proposition}
%\begin{proof} By \eqref{rel_morse<} and \eqref{rel_morse}, it suffices to show that ${\dim }(V^-\cap H^+)\ge \mlin(u,v)=\dim W^-$. If not, there is $(\phi,\psi)\in W^- \cap (V^-\cap H^+)^{\perp}$, $(\phi,\psi)\neq (0,0)$. By Proposition \ref{prop:solosym+} $\phi=\psi$, hence $\Qlin(\phi,\psi)=\Qi(\phi,\psi)$ thanks to \eqref{Qlin=Qi}, and obviously $\Qlin(\phi,\psi)<0$ as $(\phi,\psi)\in W^-$. \\ Come potrei far vedere che $\Qi(\phi,\psi)\ge 0 ?$ In generale NON \`e vero che $H^+\cap (V^-\cap H^+)^\perp \subset H^+ \cap V^+$. Il secondo insieme pu\`o essere vuoto.\\ Immaginiamo che, sotto un'ipotesi tipo "nondegenerazione" si possa dimostrare che there is $(\phi,\phi)\in W^-\setminus (0,0)$ with $0=\langle\mathcal I''(\phi,\phi), (\psi, \psi)\rangle =\Bsym(\phi,\phi), (\psi, \psi)$ per ogni $(\psi,\psi ) \in V^-\cap H^+$. \end{proof}}

% Next, let $H^+$ and $H^-$ the half-spaces \[ H^{\pm}:=\left\{ (\phi,\pm\phi) \, : \, \phi\in H^1_0 \right\}.\]

\section{Radial solutions: a decomposition formula for the Morse index}

In the present section we take that $\Omega$ is radially symmetric (a ball or a spherical cell) %$= \{ x\in \R^N : |x|<R\}$ is a ball centered at the origin 
and $H=H(|x|, u, v)$.  If $(u,v)$ is a radial solution, one can look at the restriction of the quadratic form $\Qlin$ to the space of radial functions and define the \emph{radial linearized Morse index} $\mlinrad$  as the maximal dimension of a subspace of $\left(H^1_{0,\rad}\right)^2$ where $\Qlin$ is negative, and the \emph{radial symmetric eigenvalues} $\L_{\rad}^{\sym}$ by minimizing the Rayleigh quotient first on $\left(H^1_{0,\rad}\right)^2$ and then on its subspaces. The eigenvalues are characterized by means of a radial differential problem
\begin{equation}\label{PA_sym_rad}
\begin{cases} 
-\left( r^{N-1} \phi'\right)'- r^{N-1} \left(H_{uv}  \phi +
\frac{1}{2}\Delta H \psi\right)= r^{N-1} \Lambda  \phi , \\
-\left( r^{N-1} \psi'\right)'- r^{N-1} \left( 
\frac{1}{2}\Delta H \phi + H_{uv}  \psi \right)= r^{N-1} \Lambda  \psi  , \\
\phi , \psi \in H^1_{0,\rad} , \end{cases}
\end{equation}
that can be projected onto $H^{\pm}_{\rad}$ giving rise to the decoupled system
\[\begin{cases}
-\left( r^{N-1} \xi'\right)'+ r^{N-1} a \xi= r^{N-1}\Lambda \xi , \\
-\left( r^{N-1} \eta'\right)'+ r^{N-1} b \eta= r^{N-1}\Lambda \eta  ,	\\
	\xi , \eta\in H^1_{0,\rad} 
\end{cases} \]
with $a$ and $b$ as in \eqref{def_a}, \eqref{def_b}.
In particular, the analogous of Proposition \ref{prop:solosym+} holds.
\begin{proposition}\label{prop:solosym+rad}
Let the set $ \Omega$ be radially symmetric, $H=H(|x|, u, v)$, and take $(u,v)$ a continuous radial solution of \eqref{Ham}. Under assumption \eqref{H1}, the radial linearized Morse index is equal to the number of negative eigenvalues of 
	\begin{equation}\label{PA_sym+rad}\begin{cases}
	-\left( r^{N-1} \xi'\right)'+ r^{N-1} a \xi= r^{N-1}\mu \xi  , \\
\xi \in H^1_{0,\rad} , 
	\end{cases}\end{equation}
	each counted with multiplicity. Moreover the negative radial symmetric eigenvalues coincide with the negative eigenvalues of \eqref{PA_sym+rad}, and the related eigenfunction are $(\xi ,\xi)$, where  $\xi\in H^1_{0,\rad}$ solves \eqref{PA_sym+rad} with $\mu=\Lambda^{\rad}_{\sym}$. 
\end{proposition}
We mention in passing that \eqref{PA_sym+rad} is a Sturm-Liouville problem, hence Picone's comparison Principle applies and gives that the negative radial eigenvalues are simple, and if the $n^{th}$ radial  eigenvalue is negative, then its eigenfunctions have exactly $n$ nodal zones.

\subsection{Singular symmetric eigenvalues}

The notion of singular eigenvalues, an effective  tool introduced to deal with scalar equations, also applies to systems. Let us briefly recall the definitions and main properties, and refer to \cite{AGNARWA} and the references therein for a detailed account.
Let 
\[ \widehat{L}= \left\{ \phi\in L^2 \, : \, \int_{\Omega} \frac{1}{|x|^2} \phi^2dx <\infty \right\},\] 
endowed with the scalar product 
\[ \langle\phi,\eta\rangle_{\widehat L}:= \int_{\Omega} \frac{1}{|x|^2} \phi\eta dx \]
and the induced norm
\[ \|\phi\|_{\widehat L}  := \left(\langle\phi,\phi\rangle_{\widehat L}\right)^{\frac{1}{2}} .\]
Next, let $\widehat{H}^1_0:= H^1_0 \cap \widehat L$ and $\widehat{H}^1_{0,\rad}:= H^1_{0,\rad} \cap \widehat L$.

A number $\widehat\Lambda$ is said a \emph{singular symmetric eigenvalue} if the equation
\begin{equation}\label{PA_sym_sing}
\begin{cases} 
-\Delta \phi - H_{uv}  \phi- 
\frac{1}{2}\Delta H \psi= \frac{1}{|x|^2}\L \phi \qquad & \text{ in } \Omega  ,  \\
-\Delta  \psi -\frac{1}{2}\Delta H  \phi - H_{uv} \psi= \frac{1}{|x|^2}\L \psi& \text{ in } \Omega , 
\\
\phi , \psi \in  \widehat{H}^1_0 & 
\end{cases}\end{equation}
has a nontrivial weak solution, i.e. $(\phi,\psi) \in \left(\widehat{H}^1_0\right)^2\setminus\{(0,0)\}$ with
\[ \Bsym\left((\phi,\psi), (\xi,\eta)\right) = \widehat \Lambda \left(\langle\phi,\eta\rangle_{\widehat L} +  \langle\psi,\xi\rangle_{\widehat L} \right)\]
for every $\xi, \eta  \in \widehat{H}^1_0$. 
In that case, any function $(\phi,\psi)$ solving \eqref{PA_sym_sing} is called a singular eigenfunction related to $\widehat \Lambda$.

The same arguments of the scalar case (see, for instance, \cite[Propositions 3.1, 3.2  and Lemma 3.3]{AGNARWA}) show that if 
\begin{equation*}\label{singularRayleigh1}
\widehat{\Lambda}_1:=\inf \left\{ \dfrac{\Qlin(\phi,\psi)}{\|\phi\|_{\widehat L}^2 + \|\psi\|_{\widehat L}^2  }\, : \, \phi,\psi \in \widehat{H}^1_0 \right\}< \left(\frac{N-2}{2}\right)^2, 
\end{equation*}
then it is attained by a nontrivial function $(\phi,\psi)$ which solves \eqref{PA_sym_sing} with $\L=\widehat{\Lambda}_1$, so that $\widehat{\Lambda}_1$ and $(\phi,\psi)$ are respectively an eigenvalue and a related eigenfunction.
Next we can define iteratively
\begin{equation*}\label{singularRayleighn+1} \widehat{\Lambda}_{n}:= \inf \left\{ \dfrac{\Qlin(\phi,\psi)}{\|\phi\|_{\widehat L}^2 + \|\psi\|_{\widehat L}^2}
\, : \, \phi,\psi \in \widehat{H}^1_0 , \, 	\langle\phi,\phi_k\rangle_{\widehat L}+ \langle\psi,\psi_k\rangle_{\widehat L} =0 \text{ for } k=1,\dots n-1 \right\}  . \end{equation*}
As far as $\widehat{\Lambda}_{n} < \left(\frac{N-2}{2}\right)^2$, it is attained by a nontrivial function $(\phi,\psi )$ which solves \eqref{PA_sym_sing} and then $\widehat{\Lambda}_{n}$ is a singular eigenvalue according to the previous definition.
\\
Conversely, if $\left(\Lambda, \phi,\psi \right)\in \left(-\infty, \left(\frac{N-2}{2}\right)^2\right)\times\left(\widehat{H}^1_0\right)^2\setminus\{(0,0)\}$ solves \eqref{PA_sym_sing}, then there exists $n$ such that $\Lambda = \widehat{\Lambda}_{n}$ according to \eqref{singularRayleighn+1}.
%\\ The same arguments used for the standard eigenvalues yield that the first singular eigenfunction is simple and has positive component, if the system is fully coupled (Property 1). Eigenfunctions related to different singular eigenvalues are orthogonal in $\widehat{\mathbb L}^2$, because they solve two equations of type \eqref{singulareigenvalueproblem} for different $\widehat{\Lambda}$'s (Property 2).

%An important property of singular eigenvalues is the following.\begin{proposition}\label{morseconautovalorisingolari}	The linearized Morse index is the number of negative singular symmetric eigenvalues, each counted with multiplicity.\end{proposition} The proof is exactly the same of \cite{GGN2}, and uses the continuity of the eigenvalues with respect to the set $\Omega$.

When the set $\Omega$ is radially symmetric, $H=H(|x|, u, v)$, and $(u , v)$ is a radial solution,  one can define the singular radial symmetric eigenvalues and eigenfunctions as solutions to the problem
\begin{equation}\label{PA_sym_sing_rad}
\begin{cases} 
-\left( r^{N-1} \phi'\right)'- r^{N-1} \left(H_{uv}  \phi +
\frac{1}{2}\Delta H \psi\right)= r^{N-3} \widehat{\Lambda}^{\rad}  \phi , \\
-\left( r^{N-1} \psi'\right)'- r^{N-1} \left( 
\frac{1}{2}\Delta H \phi + H_{uv}  \psi \right)= r^{N-3} \widehat{\Lambda}^{\rad} \psi , \\
\phi , \psi \in \widehat{H}^1_{0,\rad} , \end{cases}
\end{equation}
or as minima of the respective Rayleigh quotients, 
provided that the minima are below the value $\left(\frac{N-2}2\right)^2$. In particular 
\begin{equation}\label{singularRayleighrad} 
\widehat{\Lambda}_{n}^{\rad}:= \min\left\{\max\limits_{(\phi,\psi) \in V_n\setminus\{0\} } \dfrac{\Qlin(\phi,\psi)}{\|\phi\|_{\widehat L}^2 + \|\psi\|_{\widehat L}^2} : \, V_n \mbox{ n-dim. subspace of } \left(\widehat{H}^1_{0,\rad} \right)^2 \right\}  . \end{equation}

An important property of singular eigenvalues is the following.
\begin{proposition}\label{morseconautovalorisingolari}
	The linearized Morse index and the radial linearized Morse index are, respectively, the number of negative singular symmetric eigenvalues and of negative singular radial symmetric eigenvalue, each counted with multiplicity.
\end{proposition}
The proof is exactly the same of \cite{GGN2}, and uses the continuity of the eigenvalues with respect to the set $\Omega$.

\subsection{Decomposition formula}
The advantage  in dealing with singular eigenvalues is that they decompose in a radial and an angular part.
To be more specific, let us introduce some notations.
We write $Y_j$ for the Spherical Harmonics, i.e. the eigenfunctions of the Laplace-Beltrami operator on the sphere $S^{N-1}$. Of course the operator $\left(-\Delta_{S^{N-1}}\right)^{-1}$ is positive compact and selfadjoint in $L^2(S^{N-1})$ and  so  it admits a sequence of eigenvalues $0=\l_1<\l_2\leq \l_j$ and eigenfunctions $Y_j(\theta)$  which form an Hilbert basis for $L^2(S^{N-1})$. 
Namely they satisfy
\begin{equation}\label{eigen-lapl-beltr}
-\Delta_{S^{N-1}}Y_j(\theta)=\l_j Y_j(\theta)\  \text{ for }\theta\in S^{N-1} .
\end{equation}
The eigenvalues $\l_j$ are given by the well known values 
\begin{equation}\label{eigen-beltrami}
\l_j:=j(N+j-2)\ \ \text{ for }j=0,1,\dots
\end{equation} 
each of which has multiplicity 
\begin{equation}\label{multiplicity-beltrami}
N_j:=\begin{cases}
1 & \text{ when }j=0,\\
\frac{(N+2j-2)(N+j-3)!}{(N-2)!j!} & \text{ when }j\geq 1 .
\end{cases}
\end{equation}
%These eigenfunctions $Y_j(\theta)$ are bounded in $L^{\infty}(S^{N-1})$ by standard regularity theory.

The same arguments of \cite[Proposition 4.1]{AGNARWA} prove the following
\begin{proposition}\label{prop-prel-4}
	Take $\Omega$ a radially symmetric domain, $H=H(|x|, u, v)$, and $(u,v)$ a radial solution to \eqref{Ham} satisfying \eqref{Huv}.
	The singular symmetric eigenvalues $\widehat \L_n<\left(\frac{N-2}2\right)^2$ can be decomposed in radial and angular part as 
	\begin{equation}\label{decomposition_eigenvalues}
	\widehat \L _n=\widehat \L_k^{\rad}+\l_j \quad \text{ for some  $k\geq 1$ and $j\ge 0$.}
	\end{equation}
	Conversely, if a singular radial symmetric eigenvalue $\widehat \L_k^{\rad}$ according to \eqref{PA_sym_sing_rad} is such that  $ \widehat \L_k^{\rad}< \left(\frac {N-2}2\right)^2 -\l_j$ for some $j\ge 0$, then $	\widehat \L _n $ given by \eqref{decomposition_eigenvalues} is a singular symmetric eigenvalue.
	
	Moreover the set of solutions to \eqref{PA_sym_sing}  with $\widehat \L = \widehat \L _n$ is spanned by functions of type 
	\begin{equation}\label{decomposition_eigenfunctions}
	(Y_j(\theta)\phi_k(r) ,Y_j(\theta)\psi_k(r)) ,
	\end{equation}
	where $Y_j$ solves \eqref{eigen-lapl-beltr} and 
	$(\phi_k ,\psi_k) $ solves \eqref{PA_sym_sing_rad}. 
\end{proposition}

Eventually, Propositions \ref{morseconautovalorisingolari} and \ref{prop-prel-4} yield an useful formula to compute the Morse index of a radial solution $(u,v)$.
If $\widehat\L^{\rad}_k$ is a radial singular eigenvalue, we write $\textsf{m}_k$ for its multiplicity, i.e. for the dimension of the solution set of \eqref{PA_sym_sing_rad}, so that
\[ \mlinrad(u,v)= \sum\limits_{k \, : \, \widehat\L^{\rad}_k <0} \textsf{m}_k .\]
Making also use of the explicit formulas \eqref{eigen-beltrami}, \eqref{multiplicity-beltrami}, one ends up with
\begin{corollary}
Take $\Omega$ a radially symmetric domain and $H=H(|x|, u, v)$; if  $(u,v)$ is a radial solution to \eqref{Ham} satisfying \eqref{Huv}, then
\begin{align}\label{morse-index-formula-0} 
\mlin(u,v) = \sum\limits_{k \, : \, \widehat\L^{\rad}_k <0} \sum\limits_{j=0}^{M_k} \textsf{m}_k N_j , 
\end{align}
where $
M_k= \min\left\{ n\in \mathbb N  : n\ge  \sqrt{\left(\frac{N-2}{2}\right)^2-\widehat\L^{\rad}_k}-\frac{N}{2} \right\} $ and $N_j$ as in \eqref{multiplicity-beltrami}.
\end{corollary}

Next, also problems \eqref{PA_sym_sing} and \eqref{PA_sym_sing_rad} transform into decoupled systems after projection  onto $H^+$ and $H^-$. 
Let $a$ and $b$  the functions introduced in \eqref{def_a}; it is easy to check the following interesting properties.

\begin{lemma}\label{lem:PAdiag_sing}
	A number $\Lambda<\left(\frac{N-2}{2}\right)^2$ is a singular symmetric eigenvalue if and only if it is an eigenvalue for at least one between 
	\begin{equation}\label{PA_sym+_sing}\begin{cases}
	-\Delta \xi + a \xi = \dfrac{1}{|x|^2}\mu \xi \qquad & \text{ in } \Omega ,  \\
	\xi=0  & \text{ on } \partial \Omega ,
	\end{cases}\end{equation}
	\begin{equation}\label{PA_sym-_sing}\begin{cases}
	-\Delta \eta +b \eta = \dfrac{1}{|x|^2}\nu \eta  \qquad & \text{ in } \Omega ,  \\
	\eta=0  & \text{ on } \partial \Omega,
	\end{cases}\end{equation}
	and the related eigenfunctions of \eqref{PA_sym_sing} are linear combinations of functions of type $(\xi ,\xi)$ and/or $(\eta ,-\eta)$, respectively.
	
	Under the convexity assumption \eqref{H1}, the linearized Morse index of any solution $(u,v)$ is equal to the number of negative eigenvalues of the singular scalar problem \eqref{PA_sym+_sing}, each counted with multiplicity. %Moreover if $(u,v)$ is a radial solution, its radial linearized Morse index is equal to the number of negative eigenvalues of the singular scalar problem 	\begin{equation}\label{PA_sym+_sing_rad}	\begin{cases} 	-\left( r^{N-1} \xi'\right)'- r^{N-1} a \xi= r^{N-3} \mu \phi \qquad & 0<r<1 , \\	\xi \in \widehat{H}^1_{0,\rad} ,& \end{cases}	\end{equation}	each counted with multiplicity. 
\end{lemma}

If $(u,v)$ is a radial solution, one can investigate the radial Morse index by means of the radial versions of \eqref{PA_sym+_sing}, \eqref{PA_sym+_sing},   that is 
\begin{equation}\label{PA_sym+_sing_rad}
\begin{cases} 
-\left( r^{N-1} \xi'\right)'- r^{N-1} a \xi= r^{N-3} \mu \xi , & \\
\xi \in \widehat{H}^1_{0,\rad} & \end{cases}
\end{equation}
\begin{equation}\label{PA_sym-_sing_rad}\begin{cases}
-\left( r^{N-1} \eta'\right)'- r^{N-1} b \eta= r^{N-3} \nu \eta   ,  & \\
\eta\in \widehat{H}^1_{0,\rad}  .& 
\end{cases}\end{equation}
Even though the latter  Sturm-Liouville problems are singular at  $r=0$, the main comparison properties still hold (see \cite[Sections 3,4]{AGNARWA}). Next lemma immediately follows.

\begin{lemma}\label{lem-H-rad-sing}
	Consider the system \eqref{Ham} with $\Omega$ radially symmetric and $H=H(|x|, u, v)$. If $(u,v)$ is a radial solution satisfying \eqref{Huv}, then for every radial singular eigenvalue  $\widehat \Lambda_k^{\rad}<\left(\frac{N-2}{2}\right)^2$ one of the following items applies:
	\begin{enumerate}[(i)]
		\item either $\widehat \Lambda_k^{\rad}$ is simple and the respective eigenfunction is $(\phi_k, \phi_k)$, where the scalar function $\phi_k$ solves \eqref{PA_sym+_sing_rad} for $\mu=\widehat\Lambda_k^{\rad}$,
			\item or $\widehat \Lambda_k^{\rad}$ is simple and the respective eigenfunction is $(\psi_k, -\psi_k)$, where the scalar function $\psi_k$ solves \eqref{PA_sym-_sing_rad} for $\nu=\widehat\Lambda_k^{\rad}$,
		\item or $\widehat \Lambda_k$ has multiplicity two and the respective eigenspace is spanned by one eigenfunction of type (i) and one of type (ii).
		\end{enumerate}
	Moreover under the convexity assumption \eqref{H1} only item (i) may happen, so that every radial singular eigenvalue is simple and  the respective eigenfunction has the form $(\phi_k,\phi_k)$, where $\phi_k$ has exactly $k$ nodal zones. I%n particular the first eigenfunction has strictly positive components. % the  first eigenfunction has strictly positive (or negative) components.
\end{lemma}

Inserting these properties into the representation formula \eqref{morse-index-formula-0} proves Proposition \ref{prop:morse-index-formula}.
%Lastly, the radial Morse index of $(u,v)$ is the number of negative eigenvalues of \eqref{PA_sym+_sing_rad} and  the Morse index formula simplifies into \[ \mlin(u,v) = \sum\limits_{k =1}^{\mlinrad} \sum\limits_{j=0}^{M_k} N_j \] where $N_j$ and $M_k$ are given, respectively, by \eqref{multiplicity-beltrami} and \eqref{morse-index-formula-00}.

\section{Strongly coupled problems}
%\begin{equation}\label{H} 
%\begin{cases}	-\Delta u = f(v) & \text{ in } \Omega , \\ -\Delta v = g(u) \qquad & \text{ in } \Omega ,\\ u=v=0 & \text{ on } \partial\Omega. \end{cases}\end{equation}with $f, g \in C^1(\R)$, under the assumptions
%\begin{align}\label{AA}\tag{H1} f(v)/v > 0 \text{ and } g(u)/u > 0 \text{ for every } u, v \neq 0 , \\ 
%\label{BB}\tag{H2} f'(v) > 0 \text{ and } g'(u)> 0 \text{ for every } u, v \neq 0  , \\
%\label{CC}\tag{H3} f'(v) > f(v)/v  \text{ and } g'(u) > g(u)/u \text{ for every } u, v \neq 0.\end{align}
%Systems of type \eqref{H} are particular case of the Hamiltomian systems \eqref{Ham}, where the Hamiltonian has the simple form $H(u,v)=F(v)+G(u)$, being $F$ and $G$ two primitives of the functions $f$ and $g$, respectively.
%Notice that assumption \eqref{BB} implies that all the hypotheses  \eqref{H-coop}, \eqref{H-compl-acc}, \eqref{H-buono} are satisfied by every nontrivial solutions, therefore Lemma \ref{lem-H-rad-sing} applies. 

In the following we take $u , v \in H^1_{0,\rad}$ a radial solution to \eqref{Ham}, and we enforce the convexity assumption \eqref{H1} by asking that 
\begin{align}\label{A}\tag{H1}
H_u/u > 0 &  \text{ and }  H_v /v > 0  , \\
\label{B}\tag{H2}
H_{uu} > \frac{H_u - H_{uv} v}{u}  > 0 &  \text{ and }  H_{vv}> \frac{H_v - H_{uv} u}{v} > 0   , 
%\\\label{C}\tag{H3} H_{uu} >  0 &  \text{ and }  H_{vv}>  0 , 
\end{align}
for every $u , v \neq 0$. 
Though Theorem \ref{teo-stima-Morse} is proved in full only for $H$ not depending on $x$, the preliminary  properties, and in particular Proposition \ref{profilo}, hold also for $H=H(|x|, u, v)$ continuous w.r.t. $x$ and satisfying \eqref{A} and \eqref{B} for almost every $x\in \Omega$. 
%Notice that \eqref{B} strengthens \eqref{H1}. {\AL CONTROLLA RELAZIONI IPOTESI} This assumptions are satisfied, for instance, by the Lane-Emden system \begin{equation}\label{LE} \begin{cases}	-\Delta u = |v|^{q-1} v& \text{ in } \Omega , \\	-\Delta v = |u|^{p-1}u \qquad & \text{ in } \Omega ,\\ u=v=0 & \text{ on } \partial\Omega, \end{cases} \end{equation} for $p , q >1$ and, more generally, for every fully coupled system for which \[H(u,v)= F(v)+G(u)\] with $F''(v) > F'(v)/v  >0$ and $G''(u) > G'(u)/u > 0$  for every $u, v \neq 0$.
We introduce some notations before entering the details.
\begin{nota}\label{notazioni}
	In the following $\Omega$ can be either a ball $\Omega=\{x\in\R^N : |x|<R\}$, either  a spherical cell $\Omega=\{x\in\R^N : \delta<|x|<R\}$.
	We will use the notations
	\begin{itemize}
\item $s_1<\dots s_{m-1} \mbox{ and  } t_1<\dots t_{n-1}$ for the internal zero points  of $u$ and $v$, respectively,
\item $s_0=t_0=0$ or $s_0=t_0=\delta$ according if $\Omega$ is a ball or a spherical shell, and $s_m=t_n=R$,
\item $\sigma_0<\sigma_1<\dots \sigma_{m-1} \mbox{ and } \tau_0< \tau_1<\dots \tau_{n-1}$
{for the maximum/minimum points in each nodal zone of $u$ and $v$, respectively.}
\end{itemize}\end{nota}
We also point out that  in radial coordinates \eqref{Ham} reads as 	\begin{equation}\label{Ham-rad} \begin{cases}	-\left(r^{N-1} u'\right)'= r^{N-1} H_v(r,u,v)  , \\	-\left(r^{N-1} v'\right)'= r^{N-1} H_u(r,u,v)  .	\end{cases}	\end{equation}
%and if $u, v$ is a classical solution we also have	\begin{equation}\label{Ham-rad-class} \begin{cases}	-u'' - \frac{N-1}{r} u' = H_v(u,v) & 0<r < 1 , \\	-v''-\frac{N-1}{r} v'=  H_u(u,v) & 0<r < 1 , \\	u'(0)=v'(0)= 0 . &	\end{cases}	\end{equation}

Let us describe the profile of any radial solution and prove Proposition \ref{profilo}.
%\begin{proposition}\label{profilo} Let $\Omega$ be radially symmetric, $H=H(|x|, u, v)$ continuous w.r.t. $x$ and satisfying \eqref{A} for almost every $x\in \Omega$, and $(u,v)$ a classical radial solution to \eqref{Ham}. Then following holds	\begin{enumerate}[A]		\item If one between $u$ and $v$ is strictly positive (or negative), then		\begin{enumerate}[({A.}i)]			\item both $u$ and $v$ are strictly positive (or negative)			\item both $u$ and $v$ have only one critical point, which coincides with the origin if $\Omega$ is a ball.		\end{enumerate}		\item If one between $u$ and $v$ has exactly $m\ge 2$ nodal zones and $m-1$ internal zeros,  then		\begin{enumerate}[({B.}i)]			\item both $u$ and $v$ have exactly $m$ nodal zones and $m-1$ internal zeros,			\item $u$ and $v$ have the same sign in their first nodal zone (and therefore in the following ones),			\item both $u$ and $v$ have exactly one critical point inside each nodal zone. 		\end{enumerate}	\end{enumerate}	Moreover the nodal zones are intertwined, meaning that	\begin{align}\label{punti-critici}	\sigma_i, \tau_i  \in ( s_i, s_{i+1}) \cap ( t_i, t_{i+1} ) \qquad \text{ for } i=0,\dots m-1	\end{align}	and, when $m\ge 3$, also	\begin{align}	\label{zone-nodali}	s_{i-1} < t_i < s_{i+1} , \qquad t_{i-1} < s_i < t_{i+1} \qquad \text{ for } i=1,\dots m-1 .	\end{align}\end{proposition}

\begin{proof}[Proof of Proposition \ref{profilo} Part  A]	
		To fix idea, we take that $v>0$.
	
If $\Omega$ is a ball, since $v$ is smooth and radial then $v'(0)=0$. Integrating the first equation in \eqref{Ham-rad} in $(0, r)$ gives
	\begin{equation}\label{la} -r^{N-1} u'(r) = \int_{\sigma_0}^r \rho^{N-1} H_v(\rho, u(\rho), v(\rho)) d\rho > 0, \end{equation}
by assumption \eqref{A}, so that $u$ is strictly decreasing, and then, in turn, positive, and $\sigma_0$ is  the unique critical point.
Repeating the same computation with reversed role for  $v$ and $u$ shows that also $v$ as an unique critical point at the origin and completes this part of the proof.
\\
If $\Omega$ is a spherical cell, let $\sigma= \inf\left\{ r\in (\delta, r) \, : \, u'(r)=0\right\}$. By continity $u'(\sigma)=0$, and the previous argument show that 
$u$ is strictly decreasing, and then, in turn, positive on $(\sigma, R)$. On the other side, when $r\in (\delta, \sigma)$, the same computations give 
\begin{equation}\label{anello}
r^{N-1} u'(r) = \int_r^{\sigma} \rho^{N-1} H_v(\rho, u(\rho), v(\rho)) d\rho > 0,
\end{equation}
so that $u$ is strictly increasing in $(\delta, \sigma)$ and therefore positive, again.
Summing up, $u(r)$ is strictly positive on $(\delta, R)$ and has an unique critical point at $\sigma_0=\sigma$.
Repeating the same computation with reversed role for  $v$ and $u$ shows that also $v$ as an unique critical point and completes the proof of part {\it A}.
\end{proof}
\begin{proof}[Proof of Proposition \ref{profilo} Part  B] \
	\newline

{\it The case $m=2$.}
	To fix idea, we take that there exists $t\in(0,1)$ such that $v(t)=0$, $v>0$ at the left of $t$ and $v<0$ at the right of $t$.
	
	If $\Omega$ is a ball, 
the relation \eqref{la} still holds for $r\in (0,t)$, %; integrating the first equation in \eqref{Ham-rad} in $(0,r)$ and recalling that $u'(0)=0$ gives	\[ -r^{N-1} u'(r) = \int_0^r \rho^{N-1} H_v(u(\rho), v(\rho)) d\rho > 0\] by assumption \eqref{A}, 
	so that $u$ is strictly decreasing in $(0,t)$, and $0$ is a local maximum point for $u$.
	Besides $u'(r)<0$ for every  $r\in (0,1)$ may not occur, otherwise $u$ is positive on $[0,1)$ and part {\it A}  forbids $v$ to have two nodal zones.
	Hence there exists $\sigma\in (t,R)$ such that $u'<0$ on $(0,\sigma)$ and $u'(\sigma)=0$.
	Taking $r\in(\sigma, R)$ and integrating again the first equation in \eqref{Ham-rad} now gives
	\[ -r^{N-1} u'(r) = \int_{\sigma}^r \rho^{N-1} H_v(\rho, u(\rho), v(\rho)) d\rho < 0 ,\]
	so that $u$ is strictly increasing on $(\sigma, R)$.
	Therefore the function $u$ has exactly one critical point $\sigma \in (t,R)$, where the minimum is achieved. 
	Moreover $u$ can not be negative in $[0,1)$ (otherwise the case {\it A} could apply), so that $u$  is positive in the first nodal zone and changes sign exactly once.
	Eventually, switching the role of $u$ and $v$ and repeating the same reasoning allows to conclude.

	If $\Omega$ is a spherical shell, we write as before $\sigma= \inf\left\{ r\in (\delta, r) \, : \, u'(r)=0\right\}$. If $\sigma\ge t$ integrating in the interval between $(\sigma, r)$, one sees that $u$ is strictly increasing on $(\sigma, R)$, and therefore negative. Because $u(\delta)=0$ and $h$ has not any critical point in $(\delta, \sigma)$, it follows that $u$ is strictly negative on $(\delta, R)$ and case {\it A} implies that the same holds for $v$, which is not true.
	Hence $\sigma < t$, and then  $u$ is strictly decreasing on $(\sigma, t)$ and increasing on $(\delta,\sigma)$. Again, case {\it A} yields that $u$ must have a further critical point in $(t,R)$, and the conclusion follows like the ball.
	\medskip

{\it The case $m\ge 3$.}	We proceed by induction on the number of nodal zones, assuming that the claim is true for solutions with less that $m$ nodal zones and deducing that it is fulfilled also by solutions with $m$ nodal zones. \\
	Now  there exist $t_0<t_1<\dots t_{m-1} <t_m=R$ such that $v>0$ on $(t_0,t_1)$ and on any interval of type  $(t_i,t_{i+1})$ with $i$ even, and $v<0$ on on any interval of type  $(t_i,t_{i+1})$ with $i$ odd, with $v(t_i)=0$ for $i=1,\dots, m$, and also $v(t_0)=0$ if $\Omega$ is  a spherical shell. We also write $\tau_i$ for a maximum (if $i$ is even) or minimum (if $i$ is odd) point of $v$ chosen in the interval  $[t_i, t_{i+1}]$. 
	
	First we take that $\Omega$ is a ball. Reasoning as in the cases $m=2$  one sees that there exists $\sigma_1> t_1$ such that $u'(r)<0$  for $r\in (0,\sigma_1)$ and $u'(\sigma_1)=0$. 
	Let us show that 
	\begin{equation}\label{claim1}
		\sigma_1 < t_2.
		\end{equation}
	If $\sigma_1 \ge t_{m-1}$, integrating the first equation in \eqref{Ham-rad} in $(\sigma_1,r)$ for any $r>\sigma_1$ gives
	\[ -r^{N-1} u'(r) = \int_{\sigma_1}^r \rho^{N-1} H_v(\rho,u(\rho),v(\rho)) d\rho < 0 \mbox{ if $m$ is even, or } >0 \mbox{ if $m$ is odd.}\]
	In the first case, $\sigma_1$ is the only one critical point of $u$ in $(0,R)$ and it is a minimum point, so that $u$ has at most  two nodal zones. 
	In the second case, $u$ is decreasing on $(0,R)$ and therefore has only one nodal zone. Anyway, we can use the induction basis after switching the role of $u$ and $v$ and obtain the contradiction that $v$ has no more than two nodal zones.
	Hence $\sigma_1 <t_{m-1}$.
	%	If $\sigma_1 \in [t_{m-2}, t_{m-1})$, integrating the first equation in \eqref{Ham-rad} in $(\sigma_1,r)$ for any $r\in(\sigma_1, t_{m-1})$ gives	\[ -r^{N-1} u'(r) = \int_{\sigma_1}^r \rho^{N-1} H_v(\rho, u(\rho), v(\rho)) d\rho > 0 \mbox{ if $m$ is even, or } <0 \mbox{ if $m$ is odd.}\]	In the first case, $u$ is decreasing on $(0, t_{m-1})$ and only two items can occur on the interval $(t_{m-1}, R)$: either $u$ is still decreasing (so that it has only one nodal zone), or it has a further critical point $\sigma_2$ which is a global minimum (so that it has at most two nodal zones).	In the second case, $\sigma_1$ is a minimum point and two other items can occur on the interval $(t_{m-1}, 1)$: either $u$ is  still increasing (so that it has at most two nodal zone), or it has a further critical point $\sigma_2$ which is a local maximum (so that it has at most three nodal zones).	In any case, switching the role of $u$ and $v$ and applying the induction basis  gives a contradiction, so that $\sigma_1<t_{m-2}$
	If $m=3$  \eqref{claim1} is already proved, otherwise a similar reasoning shows that whenever $\sigma_1\ge t_2$, then $u$ has less than $m$ nodal zones, and the induction bases yields the contradiction that also $v$ has less then $m$ nodal zones. 
	Afterward, repeating the arguments of the case $m=2$ one sees that $u$ is strictly decreasing on $(0,\sigma_1)$ and increasing on $(\sigma_1, \sigma_2)$, $\sigma_2\in (t_2, R)$ 	being a critical point. If $m=3$, then $u$ strictly decreasing on $(\sigma_2, R)$ and the proof is completed. Otherwise one sees that $\sigma_2 <t_3$ and 
	the proof can be concluded after a finite number of steps.
\\
If $\Omega$ is a spherical shell, combining the arguments used to prove the case $m=2$ and \eqref{claim1} one sees that $\sigma= \inf\left\{ r\in (\delta, r) \, : \, u'(r)=0\right\}$ satisfies $\sigma<t_1$, and then the proof follows like the ball.
 \end{proof}

Thanks to Proposition \ref{profilo}, the two components of every radial solution have the same number of nodal zones. Henceforth we shall say \textquotedblleft a solution $(u,v)$ with $m$ nodal zones\textquotedblright, meaning that both $u$ and $v$ have $m$ nodal zones.

%%%

\subsection{Proof of Theorem \ref{teo-stima-Morse}}

We begin by estimating the radial Morse index and proving \eqref{uno}.

\begin{proposition}\label{mrad>m}
Let $\Omega$ be radially symmetric, $H=H(|x|, u, v)$ continuous w.r.t. $x$ and satisfying \eqref{A}, \eqref{B} for almost every $x\in \Omega$.
If $(u,v)$ is a classical radial solution to \eqref{Ham} with $m$ nodal zones, then its radial Morse index  is at least $m$.
\end{proposition}
\begin{proof}
	Let
	\[ \begin{array}{cc}
	%h_1(r) =\begin{cases} u(r)  & \text{ in } [0,s_1] , \\ 0 & \text{ elsewhere, }  \end{cases}  \; & \;  k_1(r) =\begin{cases} v(r) & \text{ in } [0,t_1] , \\ 0 & \text{ elswhere, } \end{cases} \; & \; W_1=(h_1, k_1), \\
	u_i(r) =\begin{cases} u(r) & \text{ in } [s_{i-1},s_i] , \\
	0 & \text{ elswhere, } \end{cases}\; & \; 
	v_i(r) =\begin{cases} v(r) & \text{ in } [t_{i-1},t_i] , \\
	0 & \text{ elswhere, } 
	\end{cases} %\; &\; W_i=(h_i, k_i) 
	\end{array}\]
	for $i=1,\dots m$. Here we have used Notations \ref{notazioni}. It is clear that $u_i, v_i\in H^1_{0,\rad}$ and $(u_i,v_i)$, $(u_j,u_j)$ are orthogonal for $i\neq j$, meaning that
	\[ \int_{s_0}^{s_m} r^{N-1} \left( u_iu_j + v_iv_j\right) dr =0 . \]
	
	It suffices to check that 
	\begin{equation}\label{claim-mrad>m}
	\Qlin (u_i,v_i) < 0 \quad \mbox{ for } i=1,\dots m.
	\end{equation}
	Using $(u_i,v_i)$ as a test function in the weak formulation of \eqref{Ham-rad} gives
	\begin{align*}
	\int_{s_0}^{s_m} r^{N-1} \left( |u'_i|^2 + |v'_i|^2\right) dr = \int_{s_0}^{s_m} r^{N-1} \left[ H_v u_i + H_u v_i \right] dr, 
	\end{align*}
	hence
	\begin{align*} 
	\Qlin (u_i,v_i) = &\int_{s_0}^{s_m}\!\!r^{N-1} \left(H_v - H_{uv} u_i - H_{vv} v_i \right)  u_i dr
	+ \int_{s_0}^{s_m} \!\!r^{N-1} \left(H_u - H_{uu} u_i - H_{uv} v_i \right)  v_i dr
	.
	%\int_0^1  r^{N-1} \!\! \left( |u'_i|^2 + |v'_i|^2\right)  dr - \int_0^1 r^{N-1} \!\! \left(   H_{uv}( u_i^2 + v_i^2) + \Delta H \, u_i v_i \right) dr \\= & 
	\end{align*}
	We only estimate the first integral of the right side:
	\begin{align*} 
\int_{s_0}^{s_m}\!\!r^{N-1} \left(H_v - H_{uv} u_i - H_{vv} v_i \right)  u_i \, dr
	=
	\int_{I} \!\!r^{N-1} \left(H_v - H_{uv} u - H_{vv} v \right)  u \, dr \\
	+
	\int_{J} \!\!r^{N-1} \left(H_v - H_{uv} u\right)  u \, dr
	\end{align*}
	where $I=(s_{i-1},s_i)\cap(t_{i-1},t_i)$ and $J= (s_{i-1},s_i)\setminus(t_{i-1},t_i)$.
	Note that on the set $I$ the functions $u$ and $v$ have the same sign by Proposition \ref{profilo}, so that
	\begin{align*}
	\int_{I} \!\!r^{N-1} \left(H_v - H_{uv} u - H_{vv} v \right)  u \, dr 
	= \int_{I} \!\!r^{N-1} \left(\frac{H_v - H_{uv} u}{v} - H_{vv} \right)  u v \, dr <0
	\end{align*}
	by assumption \eqref{B}.
	On the other hand, $J \subset (t_{i-2}, t_{i-1}) \cup (t_i, t_{i+1})$. Indeed by Proposition \ref{profilo} it is known that $(s_{i},s_{i+1})\cap (t_{i}, t_{i+1})$ is not empty, therefore $s_{i}\le t_{i+1}$, and similarly $s_{i-1}\ge t_{i-2}$.
	 Hence $uv<0$ on $J$, and 
	\begin{align*}
	\int_{J} \!\!r^{N-1} \left(H_v - H_{uv} u \right)  u \, dr 
	= \int_{J} \!\!r^{N-1} \frac{H_v - H_{uv} u}{v} u v \, dr <0
	\end{align*}
	by assumption \eqref{B}.
\end{proof}

Next we give an upper bound for the singular radial eigenvalues.

Let $\xi=u'$ and $\eta=v'$;
an easy computation shows that, if $H=H(u,v)$ does not depend on $x$, then 
\begin{equation}\label{eq-derivate} \begin{cases}
-\left(r^{N-1} \xi'\right)'- r^{N-1} \left( H_{uv} \xi + H_{vv}\eta\right) = -(N-1) r^{N-3} \xi  & 0<r < 1 , \\
-\left(r^{N-1} \eta'\right)'- r^{N-1} \left( H_{uu} \xi + H_{uv}\eta\right) = -(N-1) r^{N-3} \eta & 0<r < 1 , 
\end{cases}
\end{equation}
in weak sense, that is
\begin{equation}\label{eq-derivate-weak} \begin{split}
\int_0^1  r^{N-1} \left(\xi'\phi' + \eta'\psi' \right) dr - \int_0^1 r^{N-1} \left(   H_{uv} \left( \xi\phi+\eta\psi\right) + H_{uu} \xi\psi +H_{vv} \eta\phi \right) dr \\ = -(N-1) \int_0^1 r^{N-3} \left(\xi\phi+\eta\psi \right) dr \end{split}\end{equation}
for every $\phi,\psi\in \widehat{H}^1_{0,\rad}$.

\begin{proposition}\label{stima-autovalori-m-1}
	Let assumption \eqref{B} hold, and take $(u,v)$  a classical radial solution with $m\ge 2$ nodal zones.
	Then 
	\begin{equation}\label{tre} \widehat{\Lambda}_1^{\rad} < \dots <\widehat{\Lambda}_{m-1}^{\rad}  <- (N-1) . \end{equation}
\end{proposition}
\begin{proof}
%	Assumption \eqref{B} assures that the system is fully coupled and  Lemma \ref{lem-H-rad-sing} applies, so that the singular radial eigenvalues are simple, the first eigenfunction has positive components and each component of the $i^{th}$ eigenfunction has exactly $i$ nodal zones. 
Using  Notations  \ref{notazioni}, we define the auxiliary functions
\[ \begin{array}{cc}
\xi_k(r) =\begin{cases} \xi(r) & \text{ in } [\sigma_{k-1},\sigma_k] , \\
0 & \text{ elswhere, } \end{cases}\; & \; 
\eta_k(r) =\begin{cases} \eta(r) & \text{ in } [\tau_{k-1},\tau_k] , \\
0 & \text{ elswhere, } 
\end{cases} 
\end{array}\]
for $k=1\dots m-1$.
It is not difficult to check that $\xi_k , \eta_k\in \widehat{H}^1_{0,\rad}$ (see \cite[Lemma 3.2]{AGNonlin}), in addition $(\xi_k,\eta_k)$ and $(\xi_j,\eta_j)$ are orthogonal in $(\widehat{L})^2$ if $k\neq j$.
Moreover using $(\xi_k, \eta_k)$ as a test function in \eqref{eq-derivate-weak} gives
\begin{equation}\label{aux1}
 \begin{split}
\int_{s_0}^{s_1}  r^{N-1} \left(|\xi'_k|^2 + |\eta'_k|^2\right) dr = \int_{s_0}^{s_1} r^{N-1} \left(   H_{uu} \xi\eta_k +H_{vv} \eta\xi_k + H_{uv} \left( \xi_k^2+\eta_k^2\right) \right) dr \\  -(N-1) \int_{s_0}^{s_1} r^{N-3} \left(\xi_k^2 + \eta_k^2\right) dr \end{split}\end{equation}

Let us check that
\begin{equation}\label{stima1}
\Qlin(\xi_k, \eta_k)\le - (N-1) \int_{s_0}^{s_1} r^{N-3} (\xi_k^2 + \eta_k^2) dr
\end{equation}
for every $k=1, \dots, m-1$.
Thanks to \eqref{aux1} we have
\begin{align*}
\Qlin(\xi_k, \eta_k) =  \int_{s_0}^{s_1} r^{N-1}  \left[H_{uu} (\xi-\xi_k) \eta_k  + H_{vv} \xi_k (\eta - \eta_k) \right] dr
  -(N-1) \int_{s_0}^{s_1} r^{N-3} \left(\xi_k^2 + \eta_k^2\right) dr 
\end{align*}
	We compute, as an example, 
	\begin{align*}
	\int_{s_0}^{s_1} r^{N-1} H_{vv} \xi_k (\eta - \eta_k)  dr = \int_I r^{N-1} H_{vv} u' v' dr
	\end{align*}
	for $I= (\sigma_{k-1},\sigma_k) \setminus  (\tau_{k-1},\tau_k)$.
	First we consider the case $ \sigma_{k-1} < \tau_{k-1} $ and look at the sub-interval $I_k=(\sigma_{k-1} , \tau_{k-1} )$.
	Proposition \ref{profilo} implies that $u'v'<0$ in $I_k$: 
	indeed $I_k$ is contained in the $k^{th}$ nodal zone of $u$ and also of $v$, where $u$ and $v$ have the same sign, say they are positive. In particular, both $\sigma_{k-1}$ and $\tau_{k-1}$ are maximum points, therefore
	$u'(r)<0$ for $r\in (\sigma_{k-1}, \sigma_k)$ and $v'(r)>0$ for $r\in (\tau_{k-2}, \tau_{k-1})$. 
Using also the hypothesis \eqref{B} we conclude that
\[ \int_{I_k} r^{N-1} H_{vv} u' v' dr < 0 \]
	if $\sigma_{k-1} < \tau_{k-1} $, while there is no contribution from this term if $\sigma_{k-1} \ge  \tau_{k-1}$.
		Similarly the other terms  can be computed, and  \eqref{stima1} follows.
		Remark that equality holds only if $\sigma_{k-1}=\tau_{k-1}$ and  $\sigma_k= \tau_k$.
	
	In the first instance, \eqref{stima1} yields that   $\widehat{\Lambda}_1^{\rad} \le -(N-1)$.
	But if $\widehat{\Lambda}_1^{\rad} = -(N-1)$, then  the functions $(\xi_k,\eta_k)$ should solve \eqref{PA_sym+_sing_rad}, which is not possible since they are equal to zero in an interval. 
	
	If $m\ge 3$, we look at the subspace of $\left(\widehat{ H}^1_{0,\rad}\right)^2$ generated by $(\xi_k,\eta_k)$ with $k=1,\dots m-1$, that we denote by $V$.
	Since both the components of $(\xi_k,\eta_k)$ and $(\xi_j,\eta_j)$ have contiguous support when $k\neq j$,  it is clear that $V$ has dimension $m-1$ and $\Qlin(\phi,\psi)\le - (N-1) \int_0^1 r^{N-3} (\phi^2+\psi^2) dr$ for every $(\phi,\psi)\in V$. Hence  $\widehat{\Lambda}_{m-1}^{\rad} \le -(N-1)$ thanks to \eqref{singularRayleighrad}. Again, $\widehat{\Lambda}_{m-1}^{\rad} = -(N-1)$  may not occur, otherwise one function of type $(\xi_k,\eta_k)$ should solve \eqref{PA_sym+_sing_rad}.
\end{proof}

Putting the estimate \eqref{tre}  inside the Morse index formula \eqref{quattro}, and recalling also Proposition \ref{stima1}, yields \eqref{due} and conclude the proof of Theorem \ref{teo-stima-Morse}.

\end{document}